\documentclass[11pt]{article}

\usepackage[sort]{cite}        % nicer citations
\usepackage{amssymb}          % various ams stuff
\usepackage{latexsym}         % dito
\usepackage{amsmath}          % dito
\usepackage{eucal}            % nicer caligraphic fonts
\usepackage{bbm}              % nicer bbm fonts
\usepackage{theorem}

\title{Induction of representations in deformation
  quantization}
\author{\textbf{Henrique Bursztyn}\thanks{E-mail:
    henrique@math.toronto.edu}\addtocounter{footnote}{5}
  \\[0.1cm]
  Department of Mathematics\\
  University of Toronto\\
  100 St. George Street\\
  Toronto, Ontario, M5S 3G3\\
  Canada\\[0.5cm]
  \textbf{Stefan Waldmann}\thanks{E-mail:
    Stefan.Waldmann@physik.uni-freiburg.de}
  \\[0.1cm]
  Fakult{\"a}t f{\"u}r Mathematik und Physik\\
  Albert-Ludwigs-Universit{\"a}t Freiburg\\
  Physikalisches Institut\\
  Hermann Herder Stra{\ss}e 3\\
  D 79104 Freiburg\\
  Germany}

\date{April 2005\\[0.5cm]
  \small{Contribution to the Proceedings of the Keio Workshop 2004}}

\newcommand{\cc}[1]      {\overline{{#1}}}              
\newcommand{\id}         {\operatorname{\mathsf{id}}}

\newcommand{\SP}[1]      {\left\langle{#1}\right\rangle} 
\newcommand{\ring}[1]    {\mathsf{#1}}                 
\newcommand{\Unit}       {1}                  
\newcommand{\I}          {\mathrm{i}}

\newcommand{\Bimod}[5] {\sideset{^{\scriptscriptstyle{#1}}_{\scriptscriptstyle{#2}}}{^{\scriptscriptstyle{#4}}_{\scriptscriptstyle{#5}}}{\operatorname{#3}}}

\newcommand{\EA}   {\Bimod{}{}{\mathcal{E}}{}{\mathcal{A}}}
\newcommand{\BEA}  {\Bimod{}{\mathcal{B}}{\mathcal{E}}{}{\mathcal{A}}}
\newcommand{\AHD}  {\Bimod{}{\mathcal{A}}{\mathcal{H}}{}{\mathcal{D}}}

\newcommand{\IP}[4]{{\,}_{\scriptscriptstyle{#2}\!\!}\left\langle{{#1}}\right\rangle^{\scriptscriptstyle{#3}}_{\scriptscriptstyle{#4}}}

\newcommand{\SPEA}[1]  {\IP{{#1}}{}{\mathcal{E}}{\mathcal{A}}}

\newcommand{\SPD}[1]   {\IP{{#1}}{}{}{\mathcal{D}}}

\newcommand{\SPEH}[1]  {\IP{{#1}}{}{\mathcal{E}\otimes\mathcal{H}}{\mathcal{D}}}
\newcommand{\SPA}[1]   {\IP{{#1}}{}{}{\mathcal{A}}}
\newcommand{\BSP}[1]   {\IP{{#1}}{\mathcal{B}}{}{}}

\newcommand{\A}     {\mathcal{A}} 
\newcommand{\qA}    {\boldsymbol{\mathcal{A}}} 
\newcommand{\B}     {\mathcal{B}} 

\newcommand{\tensor}[1][{}]{\mathbin{\otimes_{\scriptscriptstyle{#1}}}}
\newcommand{\tensM}[1][{}] {\mathbin{\widehat{\otimes}_{\scriptscriptstyle{#1}}}}

\newcommand{\rep}[1][{}]  {\sideset{^*}{_{#1}}{\operatorname{\textrm{-}\mathsf{rep}}}}
\newcommand{\Rep}[1][{}]  {\sideset{^*}{_{#1}}{\operatorname{\textrm{-}\mathsf{Rep}}}}

\newtheorem{lemma} {Lemma} [section]
\newtheorem{proposition} [lemma] {Proposition}

\newtheorem{theorem} [lemma] {Theorem}

\newtheorem{definition}[lemma] {Definition}

\theorembodyfont{\rm} %%% This is to typeset examples
\newtheorem{example}[lemma]{Example}
\newtheorem{remark}[lemma]{Remark}

\newenvironment{proof}{{\sc Proof:}}{{\hspace*{\fill} $\square$\\}}

\begin{document}

\maketitle

%
% a little abstract
%

\begin{abstract}
    We discuss the procedure of Rieffel induction of representations
    in the framework of formal deformation quantization of Poisson
    manifolds. We focus on the central role played by algebraic
    notions of complete positivity.
\end{abstract}

%
% Introduction and Motivation
%

\section{Introduction}

In this note we describe how various concepts and constructions in the
theory of $C^*$-algebras carry over to the purely algebraic setting of
formal deformation quantization of Poisson manifolds. Our discussion
centers around the construction of induced representations, due to
Rieffel in the framework of $C^*$-algebras \cite{rieffel:1974a}, and
its interplay with notions of complete positivity.  Although this note
is mostly expository, we highlight some aspects of the theory that we
have not made explicit before.

Deformation quantization \cite{bayen.et.al:1978a} is a procedure to
construct algebras of quantum observables associated with classical
systems. More precisely, a classical phase space is a Poisson manifold
$(M,\{\cdot,\cdot\})$ and its quantization is a formal associative
deformation $\star$, also called a \textit{star product}, of the
classical observable algebra $C^\infty(M)$ in the direction of the
Poisson bracket. Here $C^\infty(M)$ denotes the algebra of
complex-valued smooth functions on $M$ and $\star$ is a
$\mathbb{C}[[\lambda]]$-bilinear associative multiplication on
$C^\infty(M)[[\lambda]]$ given by
\begin{equation}
    \label{eq:StarProduct}
    f \star g = \sum_{r=0}^\infty \lambda^r C_r(f, g),
\end{equation}
where $C_0(f, g) = fg$, $C_1(f, g) - C_1(g, f) = \I \{f, g\}$, $1
\star f = f = f \star 1$ and all $C_r$ are bidifferential operators.
The formal parameter $\lambda$ satisfies $\cc{\lambda} = \lambda$ and
plays the role of Planck's constant $\hbar$. We require $\star$ to be
a \emph{Hermitian star product}, in the sense that $\cc{f\star g} =
\cc{g} \star \cc{f}$, so that the $\mathbb{C}[[\lambda]]$-algebra
$(C^\infty(M)[[\lambda]],\star)$ acquires a $^*$-involution given by
pointwise complex conjugation.

Other quantum mechanical concepts can be defined in deformation
quantization analogously to the usual $C^*$-algebraic approach to
quantum theory. The starting point is to regard
$\mathbb{R}[[\lambda]]$ as an ordered ring by considering
$\sum_{r=r_0}^\infty \lambda^r a_r$ to be \textit{positive} if
$a_{r_0} > 0$, where $a_{r_0}$ is the first nonzero coefficient. Then
a $\mathbb{C}[[\lambda]]$-linear functional
\begin{equation}
    \label{eq:OmegaPos}
    \omega: C^\infty(M)[[\lambda]] \longrightarrow \mathbb{C}[[\lambda]]
\end{equation}
is called \emph{positive} if $\omega(\cc{f} \star f) \ge 0$ for all $f
\in C^\infty(M)[[\lambda]]$; a \textit{state} is a positive linear
functional such that $\omega(\Unit) = 1$, and the value $\omega(f)$ is
interpreted as the \emph{expectation value} of the observable $f$ in
the state $\omega$.

To implement the idea of superposition of states, one needs a notion
of representation in deformation quantization. Given a Hermitian star
product, a \textit{representation} consists of a pre-Hilbert space
$\mathcal{H}$ over $\mathbb{C}[[\lambda]]$ (here one uses the order
structure of $\mathbb{R}[[\lambda]]$ for the definition of positive
definite $\mathbb{C}[[\lambda]]$-valued inner products) on which
$(C^\infty(M)[[\lambda]],\star)$ acts by adjointable operators. Many
physically interesting examples can be found in
\cite{bordemann.waldmann:1998a}, see \cite{waldmann:2005b} for a
recent review.

As a next step, following the theory of $C^*$-algebras, one is led to
the construction of induced representations.  Recall that if $\A$ and
$\B$ are $C^*$-algebras, the procedure of Rieffel induction consists
of constructing representations of $\B$ from representations of $\A$
with the aid of a suitable $(\B,\A)$-bimodule $\BEA$ possessing an
$\A$-valued inner product $\SPA{\cdot,\cdot}$. For each
$^*$-representation of $\A$ on a Hilbert space
$(\mathcal{H},\SP{\cdot,\cdot})$, one considers the tensor product
$\mathcal{E}\otimes_{\scriptscriptstyle{\A}}\mathcal{H}$ over $\A$ and
the natural left action of $\B$ on it.  In order to turn this tensor
product into a Hilbert space carrying a representation of $\B$, the
key point is that one can combine $\SPA{\cdot,\cdot}$ and
$\SP{\cdot,\cdot}$ to produce an inner product on
$\mathcal{E}\otimes_{\scriptscriptstyle{\A}}\mathcal{H}$ uniquely
defined by
\begin{equation}
    \label{eq:induction}
    (x\otimes\phi,y\otimes \psi) 
    \mapsto \SP{\phi,\SPA{x,y}\cdot\psi},
\end{equation}
where $x,y \in \BEA$, $\phi, \psi \in \mathcal{H}$. An important point
of this construction where specific properties of $C^*$-algebras must
come into play is showing that the inner product defined by
\eqref{eq:induction} is positive, see e.g.
\cite{raeburn.williams:1998a} for a detailed discussion.

Understanding the positivity of inner products of the form
\eqref{eq:induction} in purely algebraic versions of Rieffel induction
is the heart of this note, see also \cite{bursztyn.waldmann:2001a}
\cite{bursztyn.waldmann:2003a:pre}. We will discuss Rieffel induction
in the framework of $^*$-algebras over ordered rings in Section
\ref{sec:induc}.  Applications to deformation quantization are
presented in Section \ref{sec:defq}.  The last section contains a
brief discussion on strong Morita equivalence, a notion closely
related to algebraic Rieffel induction.

%%%%%%%%%%%%%%%%%%%%%%%%%%%%%%%%%%%%%%%%%%%%%%%%%%%%%%%%%%%%%%%%%%%%%%%%%%%%%%%
\section{The general framework of $^*$-algebras over ordered rings}
\label{sec:general}

In order to give a unified treatment of $C^*$-algebras and the
$^*$-algebras over $\mathbb{C}[[\lambda]]$ defined by Hermitian star
products, we work in the following general algebraic setting, see
\cite{bursztyn.waldmann:2003a:pre} for details: we consider
$^*$-algebras $\A$ over a ring of the form $\ring{C} = \ring{R}(\I)$;
here $\ring{R}$ is an ordered ring, like e.g. $\mathbb{R}$ or
$\mathbb{R}[[\lambda]]$, so $\ring{C}$ is a ring extension of
$\ring{R}$ by a square root of $-1$.

Along the same lines of the discussion in the introduction, we define
a $\ring{C}$-linear functional $\omega: \mathcal{A} \longrightarrow
\ring{C}$ to be \textit{positive} if $\omega(a^*a) \ge 0$ for all
$a\in \A$, which makes sense since $\ring{R} \subseteq \ring{C}$ is
ordered.  An algebra element $a \in \mathcal{A}$ is called
\emph{positive} if its expectation values are all non-negative, i.e.
$\omega(a) \ge 0$ for all positive linear functionals $\omega$. These
notions agree with the usual ones, e.g., for $C^*$-algebras, and also
make sense for Hermitian star products.  We denote the set of positive
elements by $\mathcal{A}^+$. See \cite{schmuedgen:1990a} for more
general concepts of positivity in $O^*$-algebras and
\cite{waldmann:2004a} for a comparison between them.

We now pass to representations. A \emph{pre-Hilbert space}
$\mathcal{H}$ over $\ring{C}$ is a $\ring{C}$-module with a positive
definite inner product $\SP{\cdot,\cdot}: \mathcal{H} \times
\mathcal{H} \longrightarrow \ring{C}$, i.e. $\SP{\phi, \psi} =
\cc{\SP{\psi, \phi}}$, $\SP{\phi, \phi} > 0$ for $\phi \ne 0$ and
$\SP{\cdot,\cdot}$ is linear in the second argument. The
\emph{adjointable operators} from $\mathcal{H}_1$ to $\mathcal{H}_2$
are defined as $\ring{C}$-linear maps for which adjoints exist in the
usual sense. It is easy to check that when an adjoint exists, it is is
unique.  Note that, in the case of complex Hilbert spaces, the
Hellinger-Toeplitz theorem ensures that the adjointable operators
coincide with bounded operators.  The adjointable operators
$\mathfrak{B}(\mathcal{H})$ on a pre-Hilbert space $\mathcal{H}$ form
a $^*$-algebra in the natural way, and a \emph{$^*$-representation} of
a $^*$-algebra $\mathcal{A}$ over $\ring{C}$ on $\mathcal{H}$ is a
$^*$-homomorphism $\pi: \mathcal{A} \longrightarrow
\mathfrak{B}(\mathcal{H})$.  When $\A$ is unital, we assume that
$\pi(\Unit_{\mathcal{A}}) = \id_{\mathcal{H}}$.

\begin{example}
    \label{ex:GNS}
    If $\A$ is a $^*$-algebra over $\ring{C}$, an important class of
    examples of representations is given by an algebraic version of
    the GNS construction for $C^*$-algebras. Following
    \cite{bordemann.waldmann:1998a}, for each positive linear
    functional $\omega:\A \to \ring{C}$, one forms the space
    $\mathcal{H}_{\omega}:=\A/\mathcal{J}_\omega$, where
    $\mathcal{J}_\omega$ consists of elements $a\in \A$ with
    $\omega(a^*a)=0$. The space $\mathcal{H}_\omega$ is a pre-Hilbert
    space with inner product $\SP{\psi_a,\psi_b}:=\omega(a^*b)$, where
    $\psi_a$ denotes the class of $a \in \A$ in $\mathcal{H}_\omega$;
    the \textit{GNS $^*$-representation} of $\A$ on
    $\mathcal{H}_\omega$ is defined by $\pi(a)\psi_b:=\psi_{ab}$.

    This construction in deformation quantization gives rise to
    important formal representations of Hermitian star products, such
    as the Bargmann-Fock representation of Wick star products, or the
    Schr\"odinger representation of Weyl star products on cotangent
    bundles, see \cite{bordemann.waldmann:1998a}
    \cite{bordemann.neumaier.waldmann:1999a}.
\end{example}

For a $^*$-algebra $\A$ over $\ring{C}$, we define $\rep(\mathcal{A})$
to be the category whose objects are $^*$-representations of
$\mathcal{A}$ on pre-Hilbert spaces over $\ring{C}$ and with
adjointable intertwiners as morphisms. We refer to this category as
the \emph{representation category} (or \textit{representation theory})
of $\mathcal{A}$.  In these terms, the procedure of Rieffel induction,
to be discussed in the next section, can be seen as an explicit
construction of functors between representation categories.  Functors
which establish equivalence of categories of representations will be
briefly discussed in the last section.

%%%%%%%%%%%%%%%%%%%%%%%%%%%%%%%%%%%%%%%%%%%%%%%%%%%%%%%%%%%%%%%%%%%%%%%
\section{Complete positivity and algebraic Rieffel induction}
\label{sec:induc}

In order to describe Rieffel induction in the algebraic framework of
Section \ref{sec:general}, we need to consider algebraic analogs of
Hilbert $C^*$-modules, see e.g.  \cite{lance:1995a}. The reader may
consult \cite{bursztyn.waldmann:2003a:pre} for details.

Let $\A$ be a $^*$-algebra over $\ring{C}$, and let $\mathcal{E}$ be a
(right) $\A$-module (we may write $\EA$ to stress the $\A$-action). An
\textit{$\A$-valued inner product} on $\mathcal{E}$ is a
$\ring{C}$-sesquilinear map (linear in the second argument)
\begin{equation}
    \label{eq:SPHD}
    \SPA{\cdot, \cdot}: 
    \mathcal{E} \times \mathcal{E} \longrightarrow \mathcal{A},
\end{equation}
such that $\SPA{x,y} = \SP{y, x}_{\scriptscriptstyle{\mathcal{A}}}^*$
and $\SPA{x, y \cdot a} = \SPA{x, y} a$ for all $x, y \in \mathcal{E}$
and $a\in \mathcal{A}$. We call $\SPA{\cdot, \cdot}$
\emph{non-degenerate} if $\SPA{x, y} = 0$ for all $x$ implies $y = 0$,
in which case the pair $(\mathcal{E},\SPA{\cdot, \cdot})$ is called an
\textit{inner-product $\A$-module}.  The inner product $\SPA{\cdot,
  \cdot}$ is called \textit{positive} if $\SPA{x,x}\in \A^+$.
Finally, $\SPA{\cdot, \cdot}$ is called \emph{strongly non-degenerate}
if the map $\mathcal{E} \ni x \mapsto \SPA{x, \cdot} \in
\mathrm{Hom}_{\A}(\mathcal{E},\A)$ is a bijection.  Similar
definitions hold for left modules (the only difference is that we have
$\ring{C}$- and $\mathcal{A}$-linearity in the \emph{first} argument).

If $\B$ is another $^*$-algebra over $\ring{C}$, then a
\textit{$(\B,\A)$-inner-product bimodule} is an inner-product
$\A$-module $(\mathcal{E},\SPA{\cdot, \cdot})$ together with a
$^*$-homomorphism $\B\to \mathfrak{B}(\mathcal{E})$, where
$\mathfrak{B}(\mathcal{E})$ is the $^*$-algebra of adjointable
operators with respect to $\SPA{\cdot, \cdot}$.  Consider an object in
$\rep(\A)$, i.e., a pre-Hilbert space $(\mathcal{H},\SP{\cdot,\cdot})$
carrying a $^*$-representation of $\A$.  In order to obtain an object
in $\rep(\B)$ from $\BEA$ and $\mathcal{H}$, we follow
\cite{rieffel:1974a} and consider the algebraic tensor product
$\mathcal{E} \tensor[\A] \mathcal{H}$, which carries a left
$\B$-action, equipped with the inner product determined by
\begin{equation}
    \label{eq:RieffelInnerProduct1}
    (x \otimes \phi, y \otimes \psi) 
    \mapsto  \SP{\phi, \SPA{x, y} \cdot \psi},
\end{equation}
for $x, y \in \mathcal{E}$ and $\phi, \psi \in \mathcal{H}$. In the
framework of $C^*$-algebras, one can prove that if $\SPA{\cdot,\cdot}$
is positive, then so is the induced inner product
\eqref{eq:RieffelInnerProduct1} (see e.g
\cite{lance:1995a,raeburn.williams:1998a}). The following proposition
indicates what is algebraically needed in general.
\begin{proposition}
    \label{proposition:PequivalentCP}
    Let us assume, for simplicity, that $\A$ and $\B$ are unital, and
    let $(\BEA,\SPA{\cdot,\cdot})$ be a $(\B,\A)$-inner-product
    bimodule. Then the following are equivalent:
    \begin{enumerate}
    \item The inner product \eqref{eq:RieffelInnerProduct1} is
        positive for any $^*$-representation of $\A$.
    \item For all $n$ and all $x_1, \ldots, x_n \in \mathcal{E}$, the
        matrix $(\SPA{x_i, x_j})$ is a positive element in $
        M_n(\mathcal{A})$ (viewing $M_n(\A)$ as a $^*$-algebra over
        $\ring{C}$ in the natural way).
    \end{enumerate}
\end{proposition}
For the proof, we need the following simple lemma:
\begin{lemma}
    \label{lemma:Omegapiphiij}
    Let $\mathcal{A}$ be unital. If $\Omega: M_n(\mathcal{A})
    \longrightarrow \ring{C}$ is a positive linear functional then
    there exists a $^*$-representation $(\mathcal{H}, \pi)$ of
    $\mathcal{A}$ and vectors $\phi_1, \ldots, \phi_n \in \mathcal{H}$
    such that
    \begin{equation}
        \label{eq:nOmega}
        n \Omega(A) = \sum_{i,j} \SP{\phi_i, \pi(a_{ij})\phi_j}
    \end{equation}
    where $A = (a_{ij}) \in M_n(\mathcal{A})$. Conversely, for any
    $^*$-representation $(\mathcal{H}, \pi)$ of $\mathcal{A}$ and any
    choice of vectors $\phi_1, \ldots, \phi_n \in \mathcal{H}$, the
    right hand side of \eqref{eq:nOmega} defines a positive linear
    functional of $M_n(\mathcal{A})$ (and this defines a positive
    $\Omega$ if $1/n \in \ring{C}$).
\end{lemma}
\begin{proof}
    This is a simple application of the GNS construction and should be
    well-known. For the reader's convenience we outline the proof. Let
    $E_{ij} \in M_n(\mathcal{A})$ be the elementary matrices with
    $\Unit$ at the $(i,j)$-position and $0$ elsewhere. Then $n A =
    \sum_{i,j,k,l} E_{ji}^* a_{il} E_{kl}$. Now let
    $(\mathcal{H}_\Omega, \Pi_\Omega)$ be the GNS representation of
    $M_n(\mathcal{A})$ with respect to $\Omega$. Then define $\phi_i =
    \sum_{j} \psi_{E_{ji}} \in \mathcal{H}_\Omega$. Clearly $\pi(a) :=
    \Pi_\Omega(a E_{11})$ is a $^*$-representation of $\mathcal{A}$ on
    $\mathcal{H}_\Omega$ and we now have
    \[
    n \Omega(A) 
    = n \SP{\psi_{\Unit_{n \times n}}, 
      \Pi_\Omega(A) \psi_{\Unit_{n \times n}}}
    = \sum_{i, j} \SP{\phi_i, \pi(a_{ij}) \phi_j}.
    \]
    The converse statement can be easily checked.
\end{proof}
\begin{proof}
    We can now complete the proof of the proposition.  Let $A =
    (\SP{x_i, x_j}) \in M_n(\mathcal{A})$. Using the assumption of
    $(1)$ and the lemma, we have $\Omega(A) \ge 0$ for all positive
    linear functionals $\Omega: M_n(\mathcal{A}) \longrightarrow
    \ring{C}$.  The converse implication will follow in much more
    generality in Theorem~\ref{theorem:CompletePos}.
\end{proof}

An $\A$-valued inner product on $\mathcal{E}$ satisfying the condition
in $(2)$ is called \emph{completely positive}
\cite{bursztyn.waldmann:2003a:pre}.  If $\SPA{\cdot,\cdot}$ is
completely positive, we call $(\mathcal{E},\SPA{\cdot, \cdot})$ a
\textit{pre-Hilbert $\A$-module}; a $(\B,\A)$-inner product bimodule
for which the $\A$-valued inner product is completely positive is
called a \textit{pre-Hilbert bimodule}.
\begin{example}
    \label{ex:cpinnerprod}
    \begin{enumerate}
    \item If $(\mathcal{H},\SP{\cdot,\cdot})$ is a pre-Hilbert space
        over $\ring{C}$, then $\SP{\cdot,\cdot}$ is automatically
        completely positive;
    \item If $\A$ is a $C^*$-algebra, then any positive $\A$-valued
        inner product is completely positive;
    \item If $\A$ is a $^*$-algebra over $\ring{C}$ and $\mathcal{E}$
        is the right projective $\A$-module $P\A^n$, where $P\in
        M_n(\A)$ is a projection, then the restriction of the natural
        $\A$-valued inner product on $\A^n$ to $\mathcal{E}$ is
        completely positive.
    \item If $\A=C^\infty(M)$, then any positive strongly
        nondegenerate $\A$-valued inner product on a finitely
        generated projective (f.g.p.)  $\A$-module is completely
        positive.
    \end{enumerate}
\end{example}
To see why $(4)$ holds, note that it follows from $(3)$ that any
f.g.p.  module over $\A$ can be equipped with a completely positive
$\A$-valued inner product, and in the case where $\A=C^\infty(M)$, any
two $\A$-valued inner products on the same f.g.p. module are
equivalent. This is because, by Serre-Swan's theorem, each f.g.p.
module $\mathcal{E}$ is given by the space of sections of a complex
vector bundle $E\to M$, and strongly non-degenerate $\A$-valued inner
products on $\mathcal{E}$ correspond to hermitian fibre metrics on
$E$.  But any two such metrics on $E$ are isometric.

With the assumption of complete positivity on inner products, it turns
out that Rieffel induction can be carried out in an even broader
setting, as we now recall.

Let $\mathcal{A}$, $\mathcal{B}$ and $\mathcal{D}$ be $^*$-algebras
over $\ring{C}$ (not necessarily unital), and let
$(\BEA,\SPEA{\cdot,\cdot})$ and $(\AHD,\SPD{\cdot,\cdot})$ be right
inner-product bimodules. Let $\BEA \tensor[\mathcal{A}] \AHD$ be the
algebraic tensor product over $\A$, seen as a $(\mathcal{B},
\mathcal{D})$-bimodule in the usual way. It carries a
$\mathcal{D}$-valued inner product, generalizing
\eqref{eq:RieffelInnerProduct1}, determined by
\begin{equation}
    \label{eq:RieffelInnerProduct}
    \SPEH{x \otimes \phi, y \otimes \psi} 
    := \SPD{\phi, \SPEA{x, y} \cdot \psi},
\end{equation}
for $x, y \in \BEA$ and $\phi, \psi \in \AHD$. The main observation is
\cite{bursztyn.waldmann:2003a:pre}:
\begin{theorem}
    \label{theorem:CompletePos}
    If $\SPEA{\cdot,\cdot}$ and $\SPD{\cdot,\cdot}$ are completely
    positive, then $\SPEH{\cdot,\cdot}$ is completely positive.
\end{theorem}

To obtain a pre-Hilbert module, we
consider the quotient
\begin{equation}
    \label{eq:FtensorhatE}
    \mathcal{E} \tensM[\mathcal{A}] \mathcal{H}
    =
    \mathcal{E} \tensor[\mathcal{A}] \mathcal{H}
    \big/
    (\mathcal{E} \tensor[\mathcal{A}] \mathcal{H})^\perp,
\end{equation}
which now carries a nondegenerate, completely positive inner product
induced by \eqref{eq:RieffelInnerProduct}. So $\mathcal{E}
\tensM[\mathcal{A}] \mathcal{H}$ is a $(\mathcal{B},
\mathcal{D})$-pre-Hilbert bimodule. In fact, the tensor product
$\tensM[\mathcal{A}]$ defines a functor
\begin{equation}
    \label{eq:tensMBFunctor}
    \tensM[\mathcal{A}]:
    \rep[\mathcal{A}](\mathcal{B}) 
    \times \rep[\mathcal{D}](\mathcal{A})
    \longrightarrow
    \rep[\mathcal{D}](\mathcal{B}),
\end{equation}
where $\rep[\mathcal{D}](\mathcal{A})$ denotes the category of
$^*$-representations of $\A$ on (right) pre-Hilbert
$\mathcal{D}$-modules (in other words, pre-Hilbert
$(\A,\mathcal{D})$-bimodules). Note that, by Example
\ref{ex:cpinnerprod}, part $(1)$, if $\mathcal{D}=\ring{C}$, then
$\rep[\mathcal{D}](\mathcal{A})$ agrees with $\rep(\A)$, the
representation category of $\A$ defined in Section \ref{sec:general}.

By fixing the bimodule $\BEA$, we obtain the \emph{Rieffel induction}
functor
\begin{equation}
    \label{eq:Rieffel}
    \mathsf{R}_{\mathcal{E}}=
    \BEA \tensM[\mathcal{A}]:
    \rep[\mathcal{D}](\mathcal{A}) 
    \longrightarrow \rep[\mathcal{D}](\mathcal{B}),
\end{equation}
which allows to compare the representation theories of $\mathcal{A}$
and $\mathcal{B}$ for any auxiliary $^*$-algebra $\mathcal{D}$.  When
$\A$, $\B$ and $\mathcal{D}$ are $C^*$-algebras, one recovers the
original construction of Rieffel after suitable topological
completions.
\begin{remark}
    Note that condition $(1)$ in Proposition
    \ref{proposition:PequivalentCP} coincides with property \textbf{P}
    used in \cite{bursztyn.waldmann:2001a} for the description of
    Rieffel induction; hence Proposition
    \ref{proposition:PequivalentCP} relates the approaches of
    \cite{bursztyn.waldmann:2003a:pre} and
    \cite{bursztyn.waldmann:2001a}.
\end{remark}
\begin{example}
    \label{ex:riefinducGNS}
    Let $\A$ be a $^*$-algebra over $C$, and let $\omega:\A\to
    \ring{C}$ be a positive linear functional. We consider $\A$ as an
    $(\A,\ring{C})$-bimodule, with $\ring{C}$-valued inner product
    $\SP{a,b}_\omega:=\omega(a^*b)$. Although this is not strictly a
    pre-Hilbert bimodule according to our definition, since
    $\SP{\cdot,\cdot}_\omega$ may be degenerate, Rieffel induction
    goes through just as well. The representation of $\A$ induced by
    the canonical representation of $\ring{C}$ on itself by left
    multiplication is the GNS representation of Example \ref{ex:GNS}.
\end{example}

%%%%%%%%%%%%%%%%%%%%%%%%%%%%%%%%%%%%%%%%%%%%%%%%%%%%%%%%%%%%%%%%%%%%%%%%%%%%%%%%%%%%%%%
\section{Rieffel induction in deformation quantization}\label{sec:defq}

In this section we discuss examples of modules over Hermitian star
products which carry completely positive inner products, and hence can
be used to implement Rieffel induction in the context of deformation
quantization.  We start by recalling how classical and quantum
positive linear functionals are related in this context.
\begin{theorem}
    \label{thm:cpdef}
    Let $\qA:=(C^\infty(M)[[\lambda]],\star)$ be a Hermitian
    deformation quantization, and let $\omega_0$ be a positive linear
    functional on $C^\infty(M)$. Then one can find $\mathbb{C}$-linear
    functionals $\omega_r:\mathcal{A}\to \mathbb{C}, \; r=1,2,\ldots,$
    so that $\omega_0+\sum_{r=1}^\infty\lambda^r\omega_r$ is a
    positive linear functional of $\boldsymbol{\mathcal{A}}$.
\end{theorem}
In other words, any Hermitian star product is a \textit{positive
  deformation} in the sense of \cite{bursztyn.waldmann:2000a}.  A
proof of this theorem can be found in
\cite{bursztyn.waldmann:2004a:pre}.

We now turn our attention to examples of pre-Hilbert modules over
Hermitian star products.  Let $\mathcal{E}$ be a f.g.p. module over a
Hermitian deformation quantization
$\qA=(C^\infty(M)[[\lambda]],\star)$, and let
\[
h:\mathcal{E}\times \mathcal{E}\to \qA, \;\; (x,y)\mapsto h(x,y)
\] 
be an $\qA$-valued inner product. Then
$\mathcal{E}_0:=\mathcal{E}/(\lambda\mathcal{E})$ is a (f.g.p.) module
over $C^\infty(M)$, and $h$ naturally induces an inner product
\[
h_0:\mathcal{E}_0 \times \mathcal{E}_0 \to C^\infty(M)
\]
by $h_0([x],[y]):=h(x,y) \mod \lambda$. We refer to $h_0$ as the
\textit{classical limit} of the inner product $h$.  The next result is
an analogue in deformation quantization of Example
\ref{ex:cpinnerprod}, part $(4)$.
\begin{theorem}
    \label{thm:qmod}
    Let $\mathcal{E}$ be a f.g.p. module over a Hermitian deformation
    quantization $\qA=(C^\infty(M)[[\lambda]],\star)$, let $h$ be a
    positive, strongly nondegenerate $\qA$-valued inner product on
    $\mathcal{E}$. Then $h$ is completely positive, and its classical
    limit $h_0$ is a Hermitian fibre metric on the vector bundle $E$
    corresponding to $\mathcal{E}_0$.
\end{theorem}
\begin{proof}
    Let $h_0$ be the classical limit of $h$. We first observe that
    $h_0$ is a positive inner product on $\mathcal{E}_0$.

    Given $x \in \mathcal{E}$, consider $h_0([x],[x])=h(x,x)\mod
    \lambda \in C^\infty(M)$, and let $\omega_0$ be a positive linear
    functional on $C^\infty(M)$. By Theorem \ref{thm:cpdef}, we can
    find a positive linear functional on $\qA$ of the form $\omega=
    \omega_0 + \sum_{r\geq 1}\lambda^r \omega_r$. Since
    $\omega(h(x,x))\geq 0$ in $\mathbb{C}[[\lambda]]$, we have that
    $\omega_0(h(x,x)\mod \lambda) \geq 0$ in $\mathbb{C}$. So
    $h_0([x],[x])\geq 0$.

    A direct computation shows that $h_0$ is strongly non-degenerate.
    Thus $(\mathcal{E}_0,h_0)$ comes from a vector bundle $E$ over $M$
    carrying a Hermitian fibre metric $h_0$, and $(\mathcal{E},h)$ is
    an example of a deformation quantization of a Hermitian vector
    bundle in the sense of \cite{bursztyn.waldmann:2000b}.  By
    \cite{bursztyn.waldmann:2003a:pre}, it follows that
    $(\mathcal{E},h)$ is isometric to an $\qA$-module as the one in
    Example \ref{ex:cpinnerprod}, part $(3)$.  Hence $h$ is completely
    positive.
\end{proof}

We note that checking that the classical limit $h_0$ is strongly
nondegenerate is sufficient to guarantee the strong nondegeneracy of
$h$.  In particular, according to \cite{bursztyn.waldmann:2000b}, any
Hermitian vector bundle over $M$ can be deformed into a pre-Hilbert
module over $\qA$ which can be used for the construction of induced
representations. In this context, line bundles over $M$ play a special
role. This is because a deformation of a line bundle $L\to M$ with
respect to a star product $\star$ defines a pre-Hilbert
\textit{bimodule} for $\star$ and another deformation quantization
$\star'$ of $M$. If $M$ is symplectic, the relationship between
$\star$ and $\star'$ is that the difference of their characteristic
classes (in the sense of e.g.~\cite{gutt.rawnsley:1999a}) is $2\pi \I
c_1(L)$ \cite{bursztyn.waldmann:2002a}, where $c_1(L)$ denotes the
first Chern class of $L$.  One can then use Rieffel induction to
transfer representations from one quantization to the other.

An interesting physical example is discussed in
\cite{bursztyn.waldmann:2002a}, where it is shown that the formal
representations of star products on cotangent bundles with a
``magnetic term'' studied in
\cite{bordemann.neumaier.pflaum.waldmann:2003a} can be obtained by
Rieffel induction of the formal Schr\"odinger representation of the
standard Weyl star product. Here, the pre-Hilbert bimodule used to
implement the induction is a deformation of the line bundle associated
with a magnetic charge satisfying Dirac's quantization condition; see
\cite{waldmann:2005b} for a detailed physical discussion of this
example.

%%%%%%%%%%%%%%%%%%%%%%%%%%%%%%%%%%%%%%%%%%%%%%%%%%%%%%%%%%%%%%%%%%%%%%%%%%%%%%%%%%%%%
\section{A unified view of strong Morita equivalence}

We now briefly recall how to obtain an equivalence of categories of
representations using the functor \eqref{eq:Rieffel}.  This leads 
to a generalization of the notion of strong Morita equivalence in
$C^*$-algebras to the algebraic framework of Section
\ref{sec:general}; details can be found in
\cite{bursztyn.waldmann:2003a:pre}.
\begin{definition}
    \label{definition:SME}
    Let $\mathcal{A}$, $\mathcal{B}$ be $^*$-algebras over $\ring{C}$
    and $\BEA$ a $(\mathcal{B}, \mathcal{A})$-bimodule so that
    $\B\cdot \mathcal{E}=\mathcal{E}$ and $\mathcal{E}\cdot \A =
    \mathcal{E}$. Suppose that $\mathcal{E}$ is equipped with
    completely positive and non-degenerate inner products
    $\SPA{\cdot,\cdot}$ and $\BSP{\cdot,\cdot}$ such that
    \begin{enumerate}
    \item $\SPA{b \cdot x, y} = \SPA{x, b^* \cdot y}$,
    \item $\BSP{x \cdot a, y} = \BSP{x, y \cdot a^*}$,
    \item $\BSP{x, y} \cdot z = x \cdot \SPA{y, z}$,
    \item $\ring{C}$-span $\{\SPA{x, y} \; | \; x, y \in \mathcal{E}\}
        = \mathcal{A}$,
    \item $\ring{C}$-span $\{\BSP{x, y} \; | \; x, y \in \mathcal{E}\}
        = \mathcal{B}$.
    \end{enumerate}
    Then $\BEA$ is called a \textit{strong Morita equivalence
      bimodule}. If there exists such a bimodule then $\mathcal{A}$
    and $\mathcal{B}$ are called \textit{strongly Morita equivalent}.
\end{definition}

As discussed in \cite{bursztyn.waldmann:2003a:pre}
\cite{bursztyn.waldmann:2001b}, one recovers Rieffel's notion of
strong Morita equivalence of $C^*$-algebras from this purely algebraic
definition by passing to minimal dense ideals.

The following theorem summarizes some of the properties of strong
Morita equivalence that have well-known counterparts in ring theory
and $C^*$-algebra theory.
\begin{theorem}
    \label{theorem:SMEequivalence}
    \begin{enumerate}
    \item Strong Morita equivalence is an equivalence relation among
        nondegenerate and idempotent $^*$-algebras over $\ring{C}$.
    \item If $\BEA$ is a strong equivalence bimodule then the Rieffel
        induction functor
        \begin{equation}
            \label{EquivalenceCat}
            \BEA \tensM[\mathcal{A}]:
            \Rep[\mathcal{D}](\mathcal{A}) \longrightarrow
            \Rep[\mathcal{D}](\mathcal{B})
        \end{equation}
        establishes an equivalence of categories for any fixed
        $^*$-algebra $\mathcal{D}$.
    \item If $\BEA$ is a strong Morita equivalence bimodule for unital
        $^*$-algebras $\mathcal{A}$ and $\mathcal{B}$, then there
        exist Hermitian dual bases $(\xi_i, \eta_i)$ and $(x_j, y_j)$,
        respectively, such that
        \begin{equation}
            \label{eq:DualBasis}
            x 
            = \sum_{i=1}^n \xi_i \cdot \SPA{\eta_i, x}
            = \sum_{j=1}^m \BSP{x, y_j} \cdot x_j
        \end{equation}
        for all $x \in \BEA$. In particular, $\BEA$ is finitely
        generated and projective as right $\mathcal{A}$-module and
        also as left $\mathcal{B}$-module.
    \item If $\A$ and $\B$ are unital, then strong Morita equivalence
        implies ring-theoretic Morita equivalence.
    \end{enumerate}
\end{theorem}

Some comments are in order. The crucial point in $(1)$ is to show
transitivity, which relies on the fact that completely positive inner
products behave well under tensor products, see Theorem
\ref{theorem:CompletePos}; in $(2)$, $\Rep[\mathcal{D}](\mathcal{A})$
denotes the subcategory of $\rep[\mathcal{D}](\mathcal{A})$ consisting
of pre-Hilbert bimodules $\AHD$ satisfying the extra nondegeneracy
condition $\A\mathcal{H}=\mathcal{H}$; property $(3)$ essentially
implies $(4)$ and is also used to show that the inner products on
equivalence bimodules of unital $^*$-algebras are strongly
nondegenerate.

In \cite{bursztyn.waldmann:2003a:pre}, we describe a class of unital
$^*$-algebras, including both unital $C^*$-algebras and Hermitan
deformation quantizations, for which the converse of part $(3)$ holds,
i.e., strong and ring-theoretic Morita equivalences define the same
equivalence relation.  The comparison between these two types of
Morita equivalence becomes more interesting at the level of Picard
group(oid)s, see \cite{bursztyn.waldmann:2003a:pre} for a discussion.

The fact that, for star products, strong Morita equivalence coincides
with Morita equivalence in the classical sense of ring theory is used
in \cite{bursztyn.waldmann:2002a} to classify strong Morita equivalent
Hermitian deformation quantizations on symplectic manifolds. An
interesting problem is to investigate the precise connection between
Morita equivalence for star products and their counterparts in
$C^*$-algebraic versions of deformation quantization.

%
% many thanks
%

\section*{Acknowledgments}

It is a pleasure to thank the organizers of the Keio workshop for the
invitation and the wonderful working atmosphere. We also thank the
participants for valuable discussions and remarks.

%
% bibliography
%

%\bibliographystyle{ewde}
%\bibliography{dqarticle,dqbook,dqprocentry,dqproceeding,preprints,misc,dqthesis}

\end{document}